   \documentclass[11pt]{article}
   \usepackage{amsmath,amsfonts,amsthm,amssymb,amscd,graphicx}
   \newcommand{\ucal}{\mathcal{U}}
   \newcommand{\vcal}{\mathcal{V}}
   \newcommand{\ccal}{\mathcal{C}}
   \newcommand{\gam}{\gamma}    \newcommand{\del}{\delta}
   \newcommand{\al}{\alpha}     \newcommand{\sig}{\sigma}
   \newcommand{\lam}{\lambda}
   \newcommand{\vareps}{\varepsilon}
   
   \newcommand{\da}{_{\alpha}}   
   \newcommand{\db}{_{\beta}}    
   \newcommand{\dgam}{_{\gamma}}   
   \newcommand{\ds}{_{\sigma}}     
   \newcommand{\dab}{_{\alpha\beta}}   
   \newcommand{\dabgam}{_{\alpha\beta\gamma}} 
   \newcommand{\dagam}{_{\alpha\gamma}}   
   \newcommand{\dbgam}{_{\beta\gamma}}   
   \newcommand{\Ob}{\operatorname{Ob}\,}
   \newcommand{\Mor}{\operatorname{Mor}\,}

   \theoremstyle{definition}
   \newtheorem{definition}{Definition}[section]
   \theoremstyle{plain}
   \newtheorem{theorem}[definition]{Theorem}

   \newcommand{\bdef}{\begin{definition}}

\begin{document}
   \title{Homotopy Transition Cocycles}
   \author{by James Wirth and Jim Stasheff}

  \maketitle


\begin{abstract}
For locally
homotopy trivial fibrations, one can define transition functions
  \[
g\dab : U\da\cap U\db \to H = H(F)
  \]
where $H$ is the monoid of homotopy equivalences of $F$ to 
itself but, instead of the cocycle
condition, one obtains only that $g\dab g\dbgam$ is homotopic to $g\dagam$ as a map of
$U\da\cap
U\db\cap U\dgam$ into $H$.  Moreover, on multiple intersections, higher homotopies arise and
are relevant to classifying the fibration.

The full theory was worked out by the first author in his
1965  Notre Dame thesis
 \cite{wirth:diss}. Here we present it using language that has been developed in the interim. We also show how this points a direction
`on beyond gerbes'.
\end{abstract}
\section{Introduction}

In the theory of fibre bundles $E\to B$, a key role is played by transition functions $g\dab
: U\da\cap U\db \to G$ with respect to an open cover $\{ U\da\}$ of $B$.  Here $G$ is the
structural group of the bundle and acts as a group of transformations on the fibre $F$.  One
of the striking properties of transition functions is the cocycle condition
  \[
g\dab g\dbgam = g\dagam \quad\mbox{on}\quad U\da\cap U\db\cap U\dgam .
  \]

For fibrations, the situation is more complicated.  Assuming the fibration is locally
homotopy trivial, one can define transition functions
  \[
g\dab : U\da\cap U\db \to H = H(F)
  \]
where $H$ is the monoid of homotopy equivalences of $F$ to itself but instead of the cocycle
condition, one obtains only that $g\dab g\dbgam$ is homotopic to $g\dagam$ as a map of
$U\da\cap
U\db\cap U\dgam$ into $H$.  Moreover, on multiple intersections, higher homotopies arise and
are relevant to classifying the fibration.

The full theory was worked out by the first author in his Notre Dame thesis
 \cite{wirth:diss}.  The
intervening years have provided a language which helps organize the technicalities, though
in no way eliminating them.  If $\{ U\da\}$ is the open covering, the disjoint union
$\coprod U\da$ can be given a rather innocuous structure of a topological category $U$,
i.e., $\Ob U = \coprod U\da$ and $\Mor U=\coprod U\da\cap U\db$
that is $x\circ y=x=y$ is defined iff $x\in U\da$, $y\in U\db$ and $x=y$.  Regarding $G$ or
$H$ as a category with one object in the standard way, the one cocycle condition says that
the transition functions define a continuous functor.  The web of higher homotopies
appropriate to a fibration are precisely equivalent to a {\em functor up to strong
homotopy}, also known as a {\em homotopy coherent functor}, which 
does arise in other contexts involving topological categories \cite{cordier-porter}.

Lest the above give the impression that we have only to translate naturally occurring
homotopies into a fancy language, we point out that some powerful topology is necessary to
construct fibrations or equivalences of fibrations from the {\em homotopy  cocycle} 
data.  In particular,
the first author's patching/glueing/ \ 
recollement (mapping cylinder) theorem, which is of fundamental importance
in  more general fibration theories, is essential \cite{wirth:diss}.

Recent developments in higher homotopy theory and especially higher
gauge theory  \cite{baez:highergauge} have inspired us
to produce this belated and somewhat updated public version
of the first author's work.  Preliminary versions and a talk at the 
University of Pennsylvania have led us to work of Breen \cite{breen:schreier,
breen:asterix} and of Simpson and Hirschowitz \cite{simpson:descente} 
which have intriguing
points of contact with Wirth's much earlier thesis. Breen was well aware at the time of \cite{breen:schreier} of the relation to higher 
homotopy theory in the context of locally homotopy trivial fibrations.  
 On the other hand, {\em gerbes} are
closely related to a special case of the homotopy transition cocycles we consider. 
We restrict our point of view to the original topological setting of 
Wirth's dissertation, leaving to the future further development of the higher homotopy cocycle point of view, especially in relation to algebraic geometry, that is,
further `pursuing stacks'.

In Section 1, we begin with a swift review of standard material about fibrations and see how
the higher order transition homotopies occur naturally.  In Section 2, we review the
realization of a topological category $\ccal$ as a space, observe that for a numerable cover
$\{ U\da\}$ of $B$, the realization $|\ucal |$ has the homotopy type of $B$ \cite{segal:classifying} and
show how a functor up to strong homotopy is sufficient to induce a map of
realizations.  
Thus a fibration $E\to B$ produces a map $B
\overset{\simeq}{\;\rightarrow}
|\ucal | \to |H| = BH$, the classifying space of $H=H(F)$ as a topological monoid.

 Our emphasis is on the cocycle point of view, although such classifying maps can
also be constructed by studying the action of the based loop
space $\Omega B$ of the base $B$ on the fibre $F$, that is,
in terms of an $A_\infty$-map of $\Omega B$ to $H(F)$ 
\cite{jds:parallel,jds:bingh}.

In Section 3, we do the topology, showing how to construct a fibration from the higher order
transition functions.  The usual universal example over $BH$ is only a quasi-fibration.
Although this could be improved to a fibration by Fuchs' technique 
\cite{fuchs:fibrations}, Wirth's construction
provides a perspicuous alternative.

In Section 4, we confront the full classification theorem:  Equivalent fibrations correspond
to homotopic functors up to strong homotopy which in turn correspond to homotopic
classifying maps.  In terms of transition functions, this appears as a direct though
complicated generalization of  cocycles up to cobounding
cocycles.

In Section 5, Wirth's concept of a ``fibration theory'' is
axiomatized.  Here too the patching theorem is crucial.

Finally, in section 6, we show the relation of our approach
to foliations and Haefliger structures. We also discuss briefly how
gerbes provide a particular
instantiation of homotopy transition cocycles of a particular
`truncated' type,
leaving for further development the relation to the work of 
Breen and of Simpson and Hirschowitz.
 
  \section{Fibrations and transition functions}

Since we wish to look at things from a homotopy invariant, not geometric,  point of view, 
a natural class of fibrations to consider is that of Dold fibrations, 
those with the WCHP (Weak Covering Homotopy Property)\cite{dold:fibrations}.

  \bdef
A map $p:E\to B$ has the WCHP if for every homotopy $H: X\times I \to B$ and
$h  :X\to E$ such that $  ph(x) = H(x,0)$, there exists a homotopy $\tilde H:
X\times [-1,1]\to E$ such that $\tilde H(x, -1) =  h(x)$ and $p\tilde H(x,t)
= H(x,t)$ for $t\in [0,1]$ while $p\tilde H(x,t) = H(x,0)$ for $t\in [-1,0].$
  \end{definition}
In other words, $H$ is covered by a homotopy whose initial position is
vertically homotopic to $h.$

On the other hand, since our emphasis will be on local data such as transition functions, it
is more appropriate to consider locally homotopy trivial fibrations.

  \bdef
A map $p: E\to B$ is locally homotopy trivial over an
open covering $\{ U\da\}$ if there exist maps $h\da$ and $k\da$ such that
the following diagrams commute:
 \begin{alignat}{3}
 p^{-1}(&U\da ) &\quad\quad\begin{matrix} h\da\\ \longrightarrow\\ \longleftarrow\\ k\da
\end{matrix}  &\quad\quad U\da\times F  \notag\\[-12pt]\\
&\searrow &&\ \ \swarrow &\quad   \notag\\[-12pt]\\
&&\quad U\da  \notag
  \end{alignat}
and $h\da$ and $k\da$ are mutual fibre homotopy inverses.
  \end{definition}

The two are related by the following:

  \begin{theorem}[Dold] 
Let $B$ be a topological space which admits a numerable covering $\{U_\alpha\}$ such
that each inclusion $U_\alpha\subset B$ is nullhomotopic, 
then $p$ has the WCHP if and only if 
$p$ is fiber homotopy trivial over each $U_\alpha$.
\end{theorem}

Let $H(F)$ be the monoid of all homotopy equivalences of $F$ to itself.   

  \bdef
The {\em transition functions} $g\dab : U\da\cap U\db \to H(F)$ are defined by the equation
  \[
h\da k\db (x,f) = \bigl( x, g\dab (x)(f)\bigr) , \qquad x\in U\da\cap U\db ,\; f\in
F.
  \]
\end{definition}
For nice $F$, $g\dab$ will be continuous with respect to the compact-open topology on
$H(F)$; otherwise we are content that the adjoint
  \[
\bar{g}\dab : U\da\cap U\db \times F\to F
  \]
is continuous.
Now consider the cocycle condition.  We have
  \[
\bigl( x, g\dab (x) g\dbgam (x)f\bigr) = h\da k\db h\db k\dgam (x,f)
  \]
which is fiber homotopic to (but not necessarily equal to) $h\da k\dgam (x,f)$.  To go any
further, we use  specific homotopies $j\da : I\times U\da \times F\da\to U\da \times F$
 from $k\da\circ h\da$ at 0 to the identity at 1.  Given them, we define
  \[
g_{\alpha\beta\gamma}: I\times U\da\cap U\db\cap U\dgam \to H(F)
  \]
by
  \[
\bigl( x, g_{\alpha\beta\gamma}(t,x)f\bigr) = h\da j\db\bigl( t, k\dgam (x,f)\bigr) .
  \]
Rather than write down the complicated formulas in general, consider four-fold
intersections.  Over $U\da\cap U\db\cap U\dgam\cap U_{\delta}$ we have the diagram of
homotopies
  \[ \begin{CD}  
 g\dab g\dbgam g_{\gam\del} @>{g\dabgam g_{\gam\del}}>> g\dagam
g_{\gam\del} \\
 @V{g\dab g_{\beta\gam\del}}VV @VV{g_{\al\gam\del}}V \\
 g\dab g_{\beta\del} @>>{g_{\al\beta\del}}> g_{\al\del}.
\end{CD}  \]
Since they are defined in terms of $j\db$ and $j\dgam$ which are independent, it is
straightforward to define 
  \[
g_{\al\beta\gam\del}: I^2\times U\da\cap U\db\cap U\dgam\cap U_{\del}\to H(F)
  \]
by
  \[
\bigl( x, g_{\al\beta\gam\del}(t,s,x)f\bigr) = h\da j\db\Bigl( t, j\dgam
\bigl( s,k_{\del}(x,f)\bigr)\Bigr) .
  \]

In general for an $n+1$-fold intersection $U\ds = U_{\al_0}\cap\cdots\cap
U_{\al_n}$ with $\sig = (\al_0,\cdots ,\al_n)$, we define
  \[
g\ds : I^{n-1}\times U\ds \to H(F)
  \]
by
  \[
\bigl( x, g\ds (t_1,\cdots ,t_{n-1}, x)f\bigr)
  = h_{\al_0}j_{\al_1}\Bigl( t_1,\cdots ,j_{\al_{n-1}}\bigl( t_{n-1},
k_{\al_n}(x_1,f)\bigr)\cdots\Bigr)
  \]
and these are compatible in the following way:  If $\sig_i = (\al_0,\cdots
,\hat{\al}_i,\cdots ,x_n)$ and $(\vareps ,i)(I^{n-1})$ is the face $t_i = \vareps
\ (\vareps = 0\mbox{ or }1),$ then for $i=1,\cdots ,n-1$:
  \begin{align*}
g\ds | (1,i)(I^{n-1})\times U\ds &= g_{\sig_i} \\
g\ds | (0,i)(I^{n-1})\times U\ds &= g_{(\al_0,\cdots ,\al_i)}g_{(\al_i,\cdots
,\al_n)}. 
  \end{align*}
\bdef 
The collection $\{g_\sigma\}$ is called a {\em homotopy transition cocycle} for\
$p:E\to B.$
  \end{definition}
These relations are the key to the categorical language we introduce next.

  \section{Topological categories and realization}
One construction of a classifying map for, e.g. G-bundles,
\cite{segal:classifying,tomDieck:klass}  
regards the covering
$\ucal = \{ U\da\}$ as a category so that the construction $B$ can be applied to give
$B\ucal$ of the homotopy type of $X$ [except for language, this was known to our
ancestors] and to interpret the transition functions as a functor so that they induce
$X\simeq B\ucal\to BG$.  The classifying property can be verified directly if we
choose the appropriate realization, Milnor's construction \cite{milnor:BG} 
of a classifying space for
a topological group, which has built in a nice ``universal'' open cover.
\vskip1ex
  \bdef  A category is a {\em topological category} if the objects and morphisms
form topological spaces such that the source, target and composition maps are continuous.
  \end{definition}
\vskip1ex
  \bdef The {\em realization} $|\ccal |$ of a topological category $\ccal$ is the following
space:  Let $\ccal_0 = \Ob \ccal , \ \ccal_1 = \Mor \ccal ,\ 
\ccal_p\subset (\Mor \ccal )^p$ consists of all $p$-tuples $(f_1,\cdots
,f_p)$ such that $f_1\circ \cdots \circ f_p$ is defined.  Consider the subset
$B\ccal\subset\Delta^{\infty}\times (\ccal_1)^{\infty}$ consisting of pairs
$(\vec{t}, \{ g_{ij}\} )$ such that

1) \quad $\vec{t}\in\Delta^{\infty}$

2) \quad $i,j$ runs over all pairs such that $t_it_j \neq 0$

3) \quad $g_{ij}\in\Mor \ccal$ and  $g_{ii}= Id_i$ and

4) \quad $g_{ij}g_{jk} = g_{ik}$ if $t_it_jt_k\neq 0$. \newline
(Thus Milnor's reference to ``the space of cocyles with values in $\ccal$.'')
Topologize this space by the limit of the quotient topologies of the maps
  \[
\Delta^n\times\ccal_n \to B\ccal
  \]
given by $(s_0,\cdots ,s_n,g_1,\cdots ,g_n) \to (\vec{t},\{ g_{ij}\} )$ with
$t_{k_j}=s_j$ for some $k_0<k_1<\cdots <k_n$ and $g_{k_ik_j} = g_{i+1}\cdots g_j$.

The universal cover of $B\ccal$ is given by $U_i = \{ t_i^{-1}(0,1]\}$ and the
$g_{ij}$ coordinates regarded as functions $U_i\cap U_j\to\ccal$ are universal
transition functions.  (Strictly speaking, the $U_i$ are only point-finite, but
following  Dold,
we can deform the original $t_i$ to functions $\bar{t}_i$ which
are locally finite so the associated $\bar{t}_i^{-1}(0,1]$ are also.)

Now an open covering $\{ U\da\}$ can be regarded as a category $\hat{\ucal}$ in a
rather trivial way. 
 For simplicity, well order the index set  or, perhaps more naturally, 
 following Segal \cite{segal:classifying}, let the objects be the points of 
the intersections $U\ds$ and let the morphisms be given by the inclusions.  
This is effectively
the barycentric subdivision and hence has a natural partial order.
Let $\Ob \ccal =
\coprod U\da$ and $\Mor \ccal =\coprod_{\al >\beta} U\da\cap U\db$ with $source (x\in U\da\cap U\db) = x\in U\da$ and $target (x\in U\da\cap U\db) = x\in U\db.$ 
Thus the composition is defined only for $x\in U\da\cap U\db\cap U\dgam , \al >\beta
>\gam$ and $\ccal_p = \coprod U_{\al_0} \cdots U_{\al_p} = \coprod U\ds$ for ordered
simplices $\sig$. We refer to $ |\ucal |$ as the {\em exploded} $X$.
See Figures 1-3.

\begin{figure}[ht] \hspace{1.0in}
\includegraphics[scale=0.6]{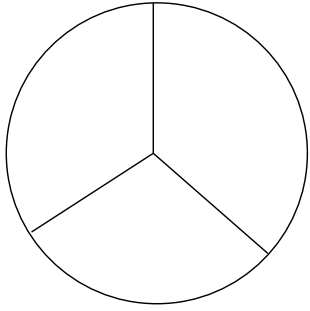}
\caption{A tri-partite covering}
\label{fig:jimfig1}
\end{figure}

\begin{figure}[ht] \hspace{1.0in}
\includegraphics[scale=0.3]{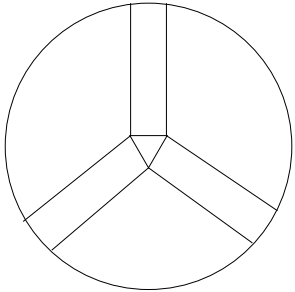}
\caption{The exploded version of Figure 1}
\label{fig:jimfig2}
\end{figure}

\begin{figure}[ht] \hspace{1.5in}
\includegraphics[scale=0.3]{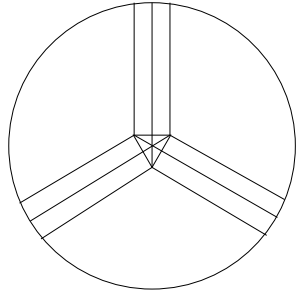}
\caption{The subdivided version of Figure 2}
\label{fig:jimfig4}
\end{figure}

There is a natural map $\pi : |\ucal |\to X$ given by forgetting the simplicial
cordinates and the multi-index.  If the covering is numerable, i.e., admits a
subordinate partition of unity $\{ p_{\al_0}(x),\cdots ,p_{\al_n}(x)\}$ for $x\in U\ds
, \sig = (\al_0,\cdots ,\al_n)$, there is a corresponding  map $\rho: X \to  |\ucal |$.  Indeed, this embeds $X$ as a deformation retract of $|\ucal |$.
See Figure 4.

\begin{figure}[ht] \hspace{1.5in}
\includegraphics[scale=0.6]{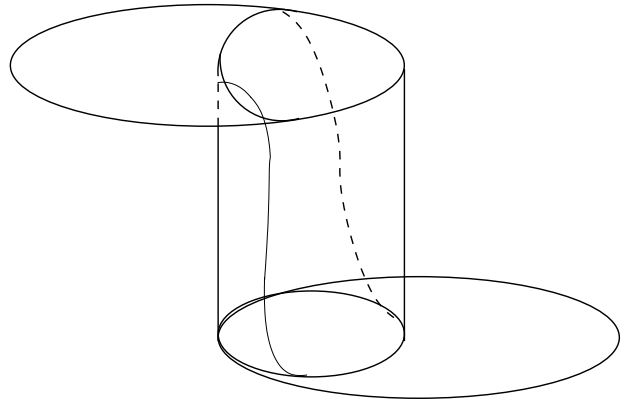}
\caption{The realization of the embedding $\rho$}
\label{fig:jimfig3}
\end{figure}
Now what is a continuous functor from $\ucal$ to $\ccal$?  It is a collection of maps
  \begin{align*}
f\da : U\da &\to \Ob \ccal \\
g\dab : U\da\cap U\db &\to \Mor \ccal
  \end{align*}
such that
  \begin{equation}
f\da = g\dab f\db
  \end{equation}
and
  \begin{equation}
 g\dab g\dbgam = g\dagam .
  \end{equation}
If $\Ob \ccal$ consists of a single point, the first condition is
trivial.  This is the case for a topological group $G = \ccal$ and gives rise to the classifying
map for a bundle
  \[
X \overset{\simeq}{\;\rightarrow} |\ucal | \to |G| = BG.
  \]
About the same time as Wirth's dissertation,
tom Dieck \cite{tomDieck:klass} proved that 
Milnor's universal
bundle $E_G$ classifies numerable $G$-bundles over arbitrary spaces
by giving explicit formulas.
Later, Wellen \cite{kapitel4} in his diplomarbeit extended the classification theorem to 
the case of groupoids $G$.

As we have seen, for a locally trivial fibration, we cannot guarantee the cocycle
condition (2) except up to homotopy.
  \end{definition}

  \bdef 
Given two topological categories $\ccal$ and $D$, a {\em functor up to strong homotopy}
(also known as  a {\em homotopy coherent functor})
is a collection of maps
  \begin{align*}
F_0 : \Ob \ccal &\to \Ob D\\
F_p : I^{p-1}\times\ccal_p &\to \Mor D
  \end{align*}
such that 
  \begin{alignat}{2}
  F_1(c: &x\to y) : F_0(x)\to F_0(y) \notag \\
F_p(t_1, &\cdots ,t_{p-1}, c_1,\cdots ,c_n) = F_{p-1}(\cdots
,\hat{t}_i,\cdots
,c_ic_{i+1},\cdots ) &&\qquad\mbox{if } t_i=0  \notag\\
&= F_i(t_1, \cdots ,t_{i-1}, c_1,\cdots ,c_i)F_{p-i}(t_{i+1},\cdots
,c_{i+1},\cdots ,c_p) &&\qquad\mbox{if } t_i=1.  \notag
  \end{alignat}
  \end{definition}
 For the special case in which the categories have one object and their
morphisms therefore form a monoid, the equivalent notion is that of a
Strongly Homotopy Multiplicative (shm) map due to Sugawara \cite{sugawara:hc} as are  the
 formulas for the corresponding map of
classifying spaces/realizations. 
The point is that $(F_0, F_1)$ does not respect the identifications on
the nose, but the higher homotopies and the connective tissue in the realization 
allow one to get around this.

  \begin{theorem}  A functor up to strong homotopy $\ccal$ to $D$
induces a map of realizations
  \[
|\ccal | \to |D|.
  \]
  \end{theorem}
In our case, our homotopy transition cocycles can be described as a 
strong homotopy functor.  Thus we have:

  \begin{theorem}
Given $p: E\to B$ locally homotopy trivial with respect to a covering
$\{ U\da\}$, a choice of coherent transition functions
  \[
g\ds : I^{n-1}\times U\ds \to H(F)
  \]
determines a map
  \[
  B \to BH(F).
  \]
  \end{theorem}

To complete the classification, we need the converse and uniqueness
up to homotopy. As for bundles, this follows from the 
construction of a universal fibration.

  \section{Construction of fibrations}

For bundles, the construction from transition functions is very easy:
$E: \coprod U\da\times F/g\dab$, i.e., for $x\in U\da\cap U\db$,
$(x,f)\in U\db\times F$ is identified with $(x,g\dab (x)f)\in
U\da\times F$.  For fibrations, not only do the transition functions
$\{ g\dab\}$ {\em not} give an equivalence relation on $\coprod
U\da\times F$ because the cocycle condition may fail, but the obvious
attempt to use the mapping cylinder $M(g\dab )$ over $U\da\cap
U\db\times I$ may only produce a quasi-fibration \cite{tulley, wirth:mapcyl}.
However Wirth has shown:

\begin{theorem} Mapping Cylinder Theorem (Wirth \cite{wirth:mapcyl})
Let $\phi : E_0\to E_1$ be a fibre homotopy equivalence over $B$,
then there is an object $\tilde{M}(\phi )$ over $I\times B$ which
serves as a mapping cylinder for $\phi$, i.e., $\tilde{M}(\phi
)|0\times B$ is $E_0$ and $\tilde{M}(\phi )|1\times B$ is $E_1$ and,
moreover, there is a {\em characterizing} homotopy equivalence $\psi$
from $\tilde{M}(\phi )$ to $I\times E_1$ such that $\psi |0\times
B=0\times\phi$ and $\psi | 1\times B=Id$. 
\end{theorem}

 In fact, the ordinary
mapping cylinder will do if $\phi$ is either a strong deformation
retract or is a fibration itself in the category of spaces over $B$.
  Finally any map over $B$ factors into a strong deformation
retract over $B$ followed by a fibration-over-$B$, just as an
ordinary map factors into a strong deformation retract followed by an
(induced from a path space) fibration.  That is, $\phi$ can be
factored as
  \[
E_0 \to E(\phi ) \to E_1
  \]
where
  \[
E(\phi ) = \{ (e,\lam )| e\in E,\lam : I\to p_1^{-1}(p_0(e)),\; \phi
(e)=\lam (0)\}.
  \]

Thus given the transition functions $\{ g\dab : U\da\times U\db \to
H(F)\}$ for $\al > \beta$, we let
  \[
E_0 = \coprod U\da\times F
  \]
and
  \[  E_1 = E_0 \cup \tilde{M}(g\dab )  \]
with the obvious identifications.  Just as we can regard $\Delta^n$
as the mapping cylinder of $\Dot{\Delta}^n\to\ast$, so we can regard
$|\ucal |$ as obtained by adding successive mapping cylinders to
$\coprod U\da$.  Let $|\ucal |(n)\subset |\ucal |$ be the subspace
represented by points with at most $(n+1)$-simplicial coordinates not
equal to zero.  Then $|\ucal |(n)=|\ucal |(n-1)\cup \coprod\ds
M(\Dot{\Delta}^n\times\ast )\times U\ds$ for all $n$-simplices
$\sig$.  Thus to extend $E_1$ to a fibration over $|\ucal |(2)$,
etc., we need to attach to $E_1$ something of the homotopy type of
$\coprod U\dabgam\times F$.  Let $h\dabgam : I\times
U\dabgam\times F\to U\dabgam\times F$ be defined by $h\dabgam
(t,x,f)=\bigl( x,g\dabgam (t,x)f\bigr)$ so that $h\dabgam$ is a fibre
homotopy equivalence over $U\dabgam$.  

Now let $E_2=E_1 \cup\coprod_{\al 
>\beta >\gam} I^2\times U\dabgam\times F$ attached by
  \begin{alignat}{2}
(s,t,x,f) &\sim (t,x,f)  &\mbox{if } s&=0 \notag\\
&\sim (s,x,f) &\mbox{if } t&=0 \notag\\
&\sim (s,x,g\dbgam (x)f) &\mbox{if } t&=1 \notag\\
&\sim (x,g\dabgam (t,x)f)\qquad&\mbox{if } s&=1. \notag
  \end{alignat}
Constructed in this way, $E_2$ will be only a quasi-fibration in general.
The proper construction is to regard the above identifications as giving a
fibre homotopy equivalence $\phi\dabgam$ of $\partial I^2\times
U\dabgam\times F$ to $E_1 | \Delta^2\times U\dabgam$ and then to attach
Wirth's $\tilde{M}(\phi\dabgam )$ instead of
 $I^2\times U\dabgam\times F$.  From Wirth's construction, there 
is a chracterizing fibre homotopy
equivalence $\psi\dabgam : I^2\times U\dabgam\times F\to
\tilde{M}(\phi\dabgam )\subset E_2$.

The general construction should now be clear.  Let $E_n$ be constructed
inductively over $|\ucal |(n)$.  Define a fibre homotopy equivalence
  \[
\psi\ds : \partial I^n\times U\ds\times F\to E_{n-1} |
\partial\Delta^n\times U\ds
  \]
by
  \begin{alignat}{2}
\psi\ds (t_1,\cdots ,t_n,x,f) &= \psi_{\sigma_i}(\cdots ,\hat{t}_i,\cdots
,x,f) &\mbox{if  } &t_i=0 \notag\\
&=\psi_{\al_0\cdots \al_{i-1}}(t_1,\cdots
,t_{i-1},x,g_{\al_{i-1}\cdots\al_n}(t_{i+1},\cdots ,t_n,x,f)) \qquad
&\mbox{if  } &t_i=1  \notag
  \end{alignat}
where $\psi\dab = g\dab$ and $\psi\da =\mbox{id}$.  Define
$E_n=E_{n-1}\cup\ds \tilde{M}(\phi\ds )$ and let
  \[
  \psi\ds : I^n\times U\ds\times F \to \tilde{M}(\phi\ds )\subset E_n
  \]
be Wirth's characterizing fibre homotopy equivalence.

In particular, this construction applies to the universal cover $\{ U_i\}$
of $BH(F)$ with transition functions
  \[
g_{ij}(\vec{t}, \{ c\dab\} ) = c_{ij}
  \]
where defined.  Notice here $g_{ij}g_{jk} = g_{jk}$ (the identity is a true
functor), but it is still important to use $\tilde{M}$ in constructing the 
universal fibration, which
is denoted $UE$ since we shall see it is the universal example
of a (WCHP) fibration with fibre $F$.

  \section{Equivalence of fibrations and transition functions}

Given that fibre homotopy equivalence is the appropriate notion, we need
to investigate the appropriate equivalence of transition functions.
Steenrod \cite{steenrod:bundles} observes that for bundles $p^i : E^i\to B, i=1,2$ with respect to two
coverings $\ucal =\{ U\da^1\}$ and $\vcal =\{ V^2\dgam\}$ with
corresponding transition functions $\{ g^1\dab\}$ and $\{
g^2_{\gam\del}\}$, the bundles are equivalent if and only if there exist maps
$\bar{g}\dagam : U\da\cap V\dgam \to G$, the group of the bundle, such
that $\{ g^1\dab ,g^2_{\gam\del}, \bar{g}\dagam\}$ satisfies the cocycle
condition.  Consider $\ucal\coprod\vcal = \{ U\da , V\dgam\}.$
The realization 
 $|\ucal\coprod\vcal |$ has $|\ucal |$ and $|\vcal |$ as deformation retracts, so
the cocycle condition above yields a homotopy
between the classifying maps between $\{ g^1\dab\}$ and $\{
g^2_{\gam\del}\}$.

Similarly we define two transition functions $\{ g\ds^i : I^{n-1}\times
U^i\ds \to H(F)\}$, $i=1,2$, to be equivalent if their union extends to a
transition function on $\ucal\coprod\vcal$.  Thus fibre homotopy
equivalent fibrations $p^i : E^i\to X$ with given fibre $F$ yield
homotopic maps $X\to BH(F)$.  [If $p^i : E^i\to X$ are locally homotopy
trivial with typical fibres $F^1 \simeq F^2$, then $BH(F^1) \simeq
BH(F^2)$ via an equivalence induced by shm maps $H(F^1)\rightleftarrows
H(F^2)$.]

Conversely, it is standard that homotopic maps induce equivalent fibrations, cf.
 one of the `Axioms of a fibration theory' 
(see section 6). 
Thus we are ready to prove:

  \begin{theorem}
For a space $F$, fibre homotopy equivalence classes of fibrations $p: E\to
X$ locally homotopy trivial with respect to numerable covers of $X$ are in
1-1 correspondence with homotopy classes of maps $X\to BH(F)$.  The
correspondence is induced by the classifying map construction above or by
the pullback of $UE$.
  \end{theorem}

We have left to show that if $f: X\to |\ucal |\to BH(F)$ is constructed
from transition functions for $p: E\to X$, then $f^* UE$ is fibre homotopy
equivalent to $E$ and that the classifying construction for $UE$ produces
a map homotopic to the identity.

Just as $|\ucal |\to X$ is induced by forgetting simplicial coordinates,
so is $f^*UE\to E$  also and, since it induces a homotopy equivalence on
each fibre, is a fibre homotopy equivalence 
\cite{dold:fibrations}.  Strictly speaking,
for $UE$ the classifying map $BH(F)\to |\ucal |\to BH(F)$ is only
homotopic to the identity since the open sets $U_i$ are
$\bar{t}_i^{-1}(0,1]$ rather than $t_i^{-1}(0,1]$.

This discussion should make it clear that a corresponding classification
theorem holds for {\em any fibration theory} as axiomatised by
Wirth \cite{wirth:diss}.

\section{Fibration theories}
A {\em fibration theory} is an assignment of a category 
${\mathcal E}  (B)$
to each topological space  $B$ and of a contravariant functor
$f^*:  {\mathcal E}(C) \to {\mathcal E}  (B)$ to each continuous map $f:B \to C$
such that $id^*$ is the identity functor and satisfying:

Axiom I: For a numerable open cover $\{U_i\}$ of a space $B$
and a system of objects (or morphisms) $\{E_i\}$ over each 
$U_i$ such that $E_i$ and $E_j$ agree over $U_i\cap U_j$,
then there exists a unique common extension of the $E_i$ over
$B$.

Axiom II: If $\phi$ is a morphism in ${\mathcal E} (B)$ such that each restriction $\phi|U_i$ for a numerable open cover $\{U_i\}$ of  $B$ is a homotopy equivalence, the $\phi$ is a homotopy equivalence.

Axiom III: If $H\in {\mathcal E}  (I\times B)$, then the restrictions $H|\{t\}\times B$ are homotopy equivalent (for objects) or homotopic (for morphisms).

Axiom IV (Mapping Cylinder Axiom) If $\phi: E\to E'\in 
{\mathcal E}  (B)$ is a homotopy equivalence, then there is an
object $M(\phi) \in {\mathcal E}  (I\times B)$ which serves as a {\em mapping cylinder} for $\phi$, that is, 
$M(\phi)$ restricts to $E$ at $t=0$ and to $E'$ at $t=1$
with a charcterising homotopy equivalence $\psi_M:M(\phi)\to
I\times E'$ which restricts to $\{0\}\times \phi$, respectively
$\{1\}\times id.$

  \section{Foliations and Gerbes}
An approach similar to that for fibrations works 
for generalized foliations or Haefliger
structures.  For a topological space $X$, a Haefliger $q$-structure 
$\{ U\da , f\da , g\dab\}$ consists of an open cover $\{ U\da\}$, maps $f\da :
U\da\to R^q$ and transition functions $g\dab : U\da\cap U\db \to
\mbox{Diffeo }R^q$ such that $f\da (\al )=g\dab (x)\circ f\db (x)$ for
$x\in U\da\cap U\db$.  These satisfy Wirth's Axiom I trivially.  The other
Axioms are modified by 1) replacing homotopy equivalence by diffeomorphism
and 2) defining equivalence of Haefliger structures to mean being induced
from a structure on $X\times I$.  Thus the above method of classification
applies.
This emphasizes the central importance of cylinder
objects and their relation to equivalence in still more general structure
theories.

The essence of all that we have said is the relation of 
local to global via patching/glueing/

\noindent recollement. 
The compatabilities are sometimes referred to as {\em descent data},
 whether in our naive topological setting
or more generally for e.g. topoi.  For example, a {\em gerbe}
can be specified by descent data in terms of an open cover
$\{U_\al\}$ and a groupoid ${\mathcal G}_\al$ for each $U_\al$
with `transition' morphisms 
 \[
g\dab \in {\mathcal G}\dab 
  \]
and morphisms 
\[
c_{\al\beta\gamma} \in {\mathcal G}_{\al\beta\gamma}
\]
acting by conjugation so that  
\[
g\dab g\dbgam = 
Ad(c_{\alpha\beta\gamma})g\dagam \quad\mbox{on}\quad U\da\cap U\db\cap U\dgam 
  \]
for an inverible element $c_{\alpha\beta\gamma}.$
The role of the homotopy $g_{\alpha\beta\gamma}$
is played by conjugation with $c_{\alpha\beta\gamma}$ and such
conjugation in a connected group corresponds to a homotopy at the
classifying space level. The higher homotopies on further multiple
intersections do not appear for gerbes since the $c_{\alpha\beta\gamma}$
satisfy a strict coherence condition.  

However, as Breen suggests,  a {\em 2-gerbe} can be specified by descent data
requiring coherence  at a higher level. There's work to be done `pursueing stacks'
and perhaps more exotic objects.

\end{document}